\documentclass[12pt]{amsart}
\usepackage{amssymb, hyperref}

\newcommand{\Bgp}{{\Z^\N}}

\long\def\forget#1\forgotten{}
\newcommand{\issuenumber}{32}
\newcommand{\issuemonth}{October}
\newcommand{\issueyear}{2011}

\newcommand{\ed}{
\newpage

\section{Unsolved problems from earlier issues}

\begin{issue}Is $\binom{\Omega}{\Gamma}=\binom{\Omega}{\Tau}$?\end{issue}
\begin{issue}Is $\ufin(\cO,\Omega)=\sfin(\Gamma,\Omega)$?And if not, does $\ufin(\cO,\Gamma)$ imply
$\sfin(\Gamma,\Omega)$?\end{issue}
\stepcounter{issue}\begin{issue}Does $\sone(\Omega,\Tau)$ imply $\ufin(\Gamma,\Gamma)$?\end{issue}
\begin{issue}Is $\fp=\fp^*$? (See the definition of $\fp^*$ in that issue.)\end{issue}
\begin{issue}Does there exist (in ZFC) an uncountable set satisfying $\sfin(\cB,\cB)$?\end{issue}
\stepcounter{issue}
\begin{issue}Does $X \nin \NON(\cM)$ and $Y\nin\mathsf{D}$ imply that $X\cup Y\nin \COF(\cM)$?\end{issue}
\begin{issue}[CH]Is $\split(\Lambda,\Lambda)$ preserved under finite unions?\end{issue}
\begin{issue}Is $\cov(\cM)=\fo$? (See the definition of $\fo$ in that issue.)\end{issue}
\stepcounter{issue}
\begin{issue}Could there be a Baire metric space $M$ of weight $\aleph_1$ and a partition
$\mathcal{U}$ of $M$ into $\aleph_1$ meager sets where for each ${\mathcal U}'\subset\mathcal U$,
$\bigcup {\mathcal U}'$ has the Baire property in $M$?\end{issue}
\stepcounter{issue} 
\begin{issue}Does there exist (in ZFC) a set of reals $X$ of cardinality $\fd$ such that all
finite powers of $X$ have Menger's property $\sfin(\cO,\cO)$?\end{issue}
\begin{issue}Can a Borel non-$\sigma$-compact group be generated by a Hurewicz subspace?\end{issue}
\begin{issue}[MA]Is there $X\sbst\bbR$ of cardinality continuum, satisfying $\sone(\BO,\BG)$?\end{issue}
\begin{issue}[CH]Is there a totally imperfect $X$ satisfying $\ufin(\cO,\Gamma)$
that can be mapped continuously onto $\Cantor$?\end{issue}
\begin{issue}[CH]Is there a Hurewicz $X$ such that $X^2$ is Menger but not Hurewicz?\end{issue}
\begin{issue}Does the Pytkeev property of $C_p(X)$ imply that $X$ has Menger's property?\end{issue}
\begin{issue}Does every hereditarily Hurewicz space satisfy $\sone(\BG,\BG)$?\end{issue}
\begin{issue}[CH]Is there a Rothberger-bounded $G\le\Bgp$ such that $G^2$ is not Menger-bounded?\end{issue}
\begin{issue}Let $\cW$ be the van der Waerden ideal. Are $\cW$-ultrafilters closed under products?\end{issue}
\begin{issue}Is the $\delta$-property equivalent to the $\gamma$-property $\binom{\Omega}{\Gamma}$?\end{issue}
\stepcounter{issue}\stepcounter{issue}
\general\end{document}}

\newcommand{\x}{\times}
\newcommand{\Cantor}{{\{0,1\}^\N}}

\newcommand{\fd}{\mathfrak{d}}
\newcommand{\fp}{\mathfrak{p}}
\newcommand{\NON}{{\mathsf   {NON}}}\newcommand{\COF}{{\mathsf   {COF}}}

\newcommand{\cM}{\mathcal{M}}
\newcommand{\cov}{\mathsf{cov}}

\newcommand{\bbR}{\mathbb{R}}
\newcommand{\fo}{\mathfrak{od}}

\renewcommand{\split}{\mathsf{Split}}\newcommand{\bq}{\begin{quote}}\newcommand{\eq}{\end{quote}}
\newcommand{\cO}{\mathcal{O}}\newcommand{\cB}{\mathcal{B}}\newcommand{\BG}{\cB_\Gamma}
\newcommand{\BO}{\cB_\Omega}

\newcommand{\sone}{\mathsf{S}_1}\newcommand{\sfin}{\mathsf{S}_\mathrm{fin}}
\newcommand{\ufin}{\mathsf{U}_\mathrm{fin}} 
\newcommand{\nin}{\not\in}\newcommand{\cU}{\mathcal{U}}
\newcommand{\cW}{\mathcal{W}}

\newcommand{\N}{\mathbb{N}}\newcommand{\Z}{\mathbb{Z}}
\newcommand{\sbst}{\subseteq}
\newcommand{\by}[2]{\par\hfill\emph{#1}, #2}\newcommand{\nby}[1]{\par\hfill\emph{#1}}\newcommand{\Tau}{\mathrm{T}}
\newcommand{\CE}{\textsc{CE}}
\newtheorem{thm}{Theorem}[section]\newcommand{\bthm}{\begin{thm}} \newcommand{\ethm}{\end{thm}}
\newtheorem{prop}[thm]{Proposition}\newcommand{\bprp}{\begin{prop}} \newcommand{\eprp}{\end{prop}}
\newtheorem{fact}[thm]{Fact}\newcommand{\bfct}{\begin{fact}} \newcommand{\efct}{\end{fact}}
\newtheorem{prob}[thm]{Problem}\newcommand{\bprb}{\begin{prob}} \newcommand{\eprb}{\end{prob}}
\newtheorem{lem}[thm]{Lemma}\newcommand{\blem}{\begin{lem}} \newcommand{\elem}{\end{lem}}
\newtheorem{claim}[thm]{Claim}\newcommand{\bclm}{\begin{claim}} \newcommand{\eclm}{\end{claim}}
\newtheorem{cor}[thm]{Corollary}\newcommand{\bcor}{\begin{cor}} \newcommand{\ecor}{\end{cor}}
\newtheorem{conj}[thm]{Conjecture}\newcommand{\bcnj}{\begin{conj}} \newcommand{\ecnj}{\end{conj}}
\theoremstyle{definition}\newtheorem{defn}[thm]{Definition}\newcommand{\bdfn}{\begin{defn}} \newcommand{\edfn}{\end{defn}}
\theoremstyle{remark}\newtheorem{rem}[thm]{Remark}\newcommand{\brem}{\begin{rem}} \newcommand{\erem}{\end{rem}}
\newtheorem{cnv}[thm]{Convention}\newcommand{\bcnv}{\begin{cnv}} \newcommand{\ecnv}{\end{cnv}}
\newtheorem{exam}[thm]{Example}\newcommand{\bexm}{\begin{exam}} \newcommand{\eexm}{\end{exam}}
\newtheorem{issue}{Issue}\newcommand{\bpf}{\begin{proof}} \newcommand{\epf}{\end{proof}}
\newcommand{\be}{\begin{enumerate}}\newcommand{\ee}{\end{enumerate}}\newcommand{\bi}{\begin{itemize}}
\newcommand{\ei}{\end{itemize}}\newcommand{\itm}{\item}
\newcommand{\general}{\small\vfill\par\noindent\hrulefill\par
\noindent\textbf{Previous issues.} The previous issues of this
bulletin are available online at\\
\url{http://front.math.ucdavis.edu/search?\&t=\%22SPM+Bulletin\%22}
\\[0.1cm]
\textbf{Contributions.} Announcements, discussions, and open problems should be emailed
to \texttt{tsaban@math.biu.ac.il}\\[0.1cm]
\textbf{Subscription.}
To receive this bulletin (free) to your e-mailbox, e-mail us.
}

\newcommand{\arXivl}[4]{\subsection{#2}{#4}\par\hfill{\arx{#1}}\par\hfill\emph{#3}}
\newcommand{\arXiv}[3]{\subsection{#2}\mbox{}\par\hfill{\arx{#1}}\par\hfill\emph{#3}}

\newcommand{\arx}[1]{\url{http://arxiv.org/abs/#1}}

\title[$\mathcal{SPM}$ Bulletin \textbf{\issuenumber} (\issuemonth{} \issueyear)]{%
$\mathcal{SPM}$ Bulletin\\[0.5cm]
Issue number \issuenumber: \issuemonth{} \issueyear{} \CE{}}

\begin{document}
\maketitle


\section{Editor's note}

Dear Friends,

\bigskip

\textbf{1.} By now probably most of you know that Misha Matveev has passed away recently.
I quote Ronnie Levy's concise description, distributed via \emph{Topology News}:
\begin{quote}
Misha Matveev of George Mason University died of an apparent heart attack.
He was found in his office at approximately 1 AM on March 17.

Misha was a prolific researcher in general topology. His interests
included, but were not limited to, star-covering properties, selection
principles (such as the Rothberger and Menger properties), and monotonic
covering properties.
\end{quote}
From Maddalena (Milena) Bonanzinga, I learned that Misha was thinking then on a research topic
for his forthcoming visit to the University of Messina.
It is symbolic that such an excellent mathematician passes away while doing mathematics.

I met Misha only once in person, in the
44th annual Spring Topology and Dynamics Conference, Mississippi State University,
Mississippi State, USA, 2010, and quickly noticed his humble and kind character.
Misha injected many fresh ideas and perspectives into the field, and
I always had in mind
the hope that one day, I will collaborate with him on some of his new ideas concerning SPM.
I have recently visited the University of Messina, one of Misha's favorite places to visit.
I was fortunate to collaborate, for my first time, with Filippo Cammaroto, Milena Bonanzinga,
Bruno Antonio Pansera, and Andrei Cataliato --- all former collaborators of Misha. I was also given
the opportunity to make comments and suggestions for a nearly finished paper of Misha with Bonanzinga.
This is of some consolation for me. I would like to use this opportunity to thank my friends from
Messina for giving me this opportunity.

\bigskip

\textbf{2.} With a sharp change from bad news to good news,
I am glad to inform you that the Fourth SPM Workshop will take place on
the coming June (2012). A preliminary announcement is given below.
Please circulate this information among your friends, students, and colleagues.

\medskip

\by{Boaz Tsaban}{tsaban@math.biu.ac.il}

\hfill \texttt{http://www.cs.biu.ac.il/\~{}tsaban}

\section{IV Workshop on Coverings, Selections and Games in Topology}

Dear colleague,

\bigskip

Next year, our mutual friend Ljubi\v{s}a Ko\v{c}inac turns 65.
For this occasion, I am organizing the
\emph{IV Workshop on Coverings, Selections and Games in Topology}.
Ljubisa Kocinac initiated and started this series of conferences in
2002, Lecce, Italy.

The workshop will take place at the
Department of Mathematics,
Seconda Università di Napoli,
Caserta, Italy.

Tentative time table:  25--30 June 2012.
(Arrival: 25, work: 26--28/29, excursion: 29, departure: 30.)

\subsection{Organizing Committee}
Agata Caserta,
Giuseppe Di Maio (chair),
Dragan Djurcic.

\subsection{Scientific Committee}
Alexander V. Arhangel'skii,
Giuseppe Di Maio,
Ljubisa D.R. Kocinac,
Masami Sakai,
Marion Scheepers,
Boaz Tsaban,
C. Guido,
R. Lucchetti.

\subsection{Further information}
Each talk will last about 30 minutes.

Of course, in a period of financial cuts we do not know up to now  the
support that
we can offer to participants, the amount of registration fee, etc.

This is indeed a very preliminary report, written to circulate this
important and happy news.
We would appreciate your forwarding this information to anyone who may
be interested
in attending this conference.

On the behalf of the organizing committee, I hope to see you in Caserta.

\nby{Giuseppe Di Maio}

\section{Long announcements}

\arXivl{1106.2235}
{Constructing universally small subsets of a given packing index in Polish groups}
{Taras Banakh and Nadya Lyaskovska}
{A subset of a Polish space $X$ is called universally small if it belongs to
each ccc $\sigma$-ideal with Borel base on $X$. Under CH in each uncountable
Abelian Polish group $G$ we construct a universally small subset $A_0\subset G$
such that $|A_0\cap gA_0|=\mathfrak c$ for each $g\in G$. For each cardinal
number $\kappa\in[5,\mathfrak c^+]$ the set $A_0$ contains a universally small
subset $A$ of $G$ with sharp packing index
$\sup\{|\mathcal D|^+:\mathcal D\subset \{gA\}_{g\in G}$
is disjoint$\}$ equal to $\kappa$.}

\arXivl{1106.3127}
{Amenability and Ramsey Theory}
{Justin Tatch Moore}
{The purpose of this article is to connect the notion of the amenability of a
discrete group with a new form of structural Ramsey theory. The Ramsey
theoretic reformulation of amenability constitutes a considerable weakening of
the F\o{}lner criterion. As a by-product, it will be shown that in any non
amenable group G, there is a subset E of G such that no finitely additive
probability measure on G measures all translates of E equally.}

\arXivl{1106.4735}
{Hindman's Theorem, Ellis's Lemma, and Thompson's group $F$}
{Justin Tatch Moore}
{The purpose of this article is to formulate generalizations of Hindman's
Theorem and Ellis's Lemma for non associative groupoids. These conjectures will
be proven to be equivalent and it will be shown that they imply the amenability
of Thompson's group F. In fact the amenability of F is equivalent to a finite
form of the conjectured extension of Hindman's Theorem.}

\arXivl{1106.5116}
{A counterexample in the theory of $D$-spaces}
{Daniel T. Soukup, Paul J. Szeptycki}
{Assuming $\diamondsuit$, we construct a $T_2$ example of a hereditarily
Lindel\"of space of size $\omega_1$ which is not a $D$-space. The example has
the property that all finite powers are also Lindel\"of.}

\arXivl{1107.5383}
{Borel's Conjecture in Topological Groups}
{Fred Galvin and Marion Scheepers}
{We introduce a natural generalization of Borel's Conjecture. For each
infinite cardinal number $\kappa$, let {\sf BC}$_{\kappa}$ denote this
generalization. Then ${\sf BC}_{\aleph_0}$ is equivalent to the classical Borel
conjecture. We obtain the following consistency results:
\be
\itm If it is consistent that there is a 1-inaccessible cardinal then it is
consistent that ${\sf BC}_{\aleph_1}$ holds.
\itm If it is consistent that ${\sf BC}_{\aleph_1}$ holds, then it is
consistent that there is an inaccessible cardinal.
\itm If it is consistent that there is a 1-inaccessible cardinal with $\omega$
inaccessible cardinals above it, then $\neg{\sf BC}_{\aleph_{\omega}} \, +\,
(\forall n<\omega){\sf BC}_{\aleph_n}$ is consistent.
\itm If it is consistent that there is a 2-huge cardinal, then it is consistent
that ${\sf BC}_{\aleph_{\omega}}$ holds.
\itm If it is consistent that there is a 3-huge cardinal, then it is consistent
that ${\sf BC}_{\kappa}$ holds for a proper class of cardinals $\kappa$ of
countable cofinality.
\ee}

\arXivl{1108.2533}
{The topology of ultrafilters as subspaces of $2^\omega$}
{Andrea Medini and David Milovich}
{Using the property of being completely Baire, countable dense homogeneity and
the perfect set property we will be able, under Martin's Axiom for countable
posets, to distinguish non-principal ultrafilters on $\omega$ up to
homeomorphism. Here, we identify ultrafilters with subpaces of $2^\omega$ in
the obvious way. Using the same methods, still under Martin's Axiom for
countable posets, we will construct a non-principal ultrafilter $\cU\subseteq
2^\omega$ such that $\cU^\omega$ is countable dense homogeneous. This
consistently answers a question of Hru\v{s}\'ak and Zamora Avil\'es. Finally,
we will give some partial results about the relation of such topological
properties with the combinatorial property of being a $\mathrm{P}$-point.}

\subsection{Another note on the class of paracompact spaces whose product
with every paracompact space is paracompact}
Abstract. The paper contains the following two results:
\be
\item Let $X$ be a paracompact space and $M$ be a metric space such that $X$ can be embedded
in $M^{\aleph_1}$ in such a way that the projections of $X$ onto initial countably many coordinates
are closed. Then the product $X\x Y$ is paracompact for every paracompact space $Y$ if and
only if the first player of the $G(DC,X)$ game, introduced by Telgarsky, has a
winning strategy.
\item If $X$ is paracompact space, $Y$ is a closed image of $X$ and the first player of the
$G(DC,X)$ game has a winning strategy, then also the first player of the $G(DC,Y)$ game
has a winning strategy.
\ee
\hfill\emph{K. Alster}

\subsection{On paracompactness in the Cartesian products and the Telgarsky's game}
Let $X$ be a paracompact space and $M$ be a metric space such that $X$ can be embedded
in $M^{\aleph_1}$ in such a way that the projections $p_\alpha:X\to M^\alpha$ are closed
at every $x\in X$, and $p_\alpha^{-1}p_\alpha(x)$ is clopen for all $x\in X$.
Then the product $X\x Y$ is paracompact for every paracompact space $Y$ if and
only if the first player of the $G(DC,X)$ game, introduced by Telgarsky, has a
winning strategy.

\hfill\emph{K. Alster}

\arXivl{1109.1736}
{Elementary chains and compact spaces with a small diagonal}
{Alan Dow and Klaas Pieter Hart}
{It is a well known open problem if, in ZFC, each compact space with a small
diagonal is metrizable. We explore properties of compact spaces with a small
diagonal using elementary chains of submodels. We prove that ccc subspaces of
such spaces have countable $\pi$-weight. We generalize a result of Gruenhage
about spaces which are metrizably fibered. Finally we discover that if there is
a Luzin set of reals, then every compact space with a small diagonal will have
many points of countable character.}

\section{Short announcements}\label{RA}

\arXiv{1106.2916}{On large indecomposable Banach spaces}{Piotr Koszmider}

\arXiv{1106.2917}{A $C(K)$ Banach space which does not have the Schroeder-Bernstein property}{Piotr Koszmider}

\arXiv{1109.5281}
{Linearly Ordered Families of Baire 1 Functions}
{M\'arton Elekes}

\arXiv{1109.5283}
{Chains of Baire class 1 functions and various notions of special trees}
{M\'arton Elekes and Juris Steprans}

\arXiv{1109.5284}
{Transfinite Sequences of Continuous and Baire Class 1 Functions}
{M\'arton Elekes and Kenneth Kunen}

\ed